\documentclass{amsart}

\usepackage{indentfirst}
\usepackage{amsmath}
\usepackage{amssymb}
\usepackage[dvipdfmx]{graphicx}
\usepackage{color}
\usepackage{amsthm}
\usepackage{mathrsfs}
\usepackage{pinlabel}
\usepackage{float}
\usepackage{amsopn}

\DeclareMathOperator{\conv}{Conv}

\theoremstyle{plain}
 \newtheorem{theo}{Theorem}[section]

 \newtheorem{prob}[theo]{Problem}

\theoremstyle{definition}
 \newtheorem{defi}[theo]{Definition}

\theoremstyle{remark}
 \newtheorem{rem}[theo]{Remark}
 \newtheorem{exam}[theo]{Example}

\newcommand{\mb}[1]{\mathbf{#1}}

\newcommand{\1}{\hspace{1mm}}

\title{A note on alternating links and root polytopes}
\author{Hiroki Murakami}
\date{}

\begin{document}

\maketitle

\begin{abstract}
 In this paper, a relationship between the determinant of an alternating link and a certain polytope obtained from the link diagram is analyzed. We also show that when the  underlying graph of the link diagram is properly oriented, the number of its spanning arborescences is equal to the determinant, i.e., the value at $-1$ of the Jones polynomial, of the link.
\end{abstract}

\section{Introduction}

For a given link $L$ with diagram $D$, we construct a bipartite graph $G$ using the checkerboard coloring.
The planar dual $G^*$ of $G$ is naturally directed so that it has spanning arborescences, that is, spanning trees of $G^*$ which are directed toward a fixed root.
The number of spanning arborescences is equal to the number of hypertrees of both hypergraphs corresponding to $G$ \cite{K}.
As hypertrees do not depend on root, this provides a new proof of the known fact \cite{S} that the number of arborescences is independent of root.

Postnikov showed \cite{P} that the number of hypertrees is proportional to the volume of the root polytope corresponding to $G$.
In this paper, we make the following connection.

\begin{theo}
 Given an alternating diagram $D$ of the link L, the determinant of $L$ is equal to
 	\begin{enumerate}
	 \item the number of hypertrees in $G$,
	 \item the number of spanning aroborescences of $G^*$.
	\end{enumerate}
\end{theo}

\noindent
{\bf Organization.}
In section $2$ we recall some definitions and prove Theorem $1.1$.

\noindent
{\bf Acknowledgements.}
The author should like to express his gratitude to Tam\'as K\'alm\'an for his constant encouragement and many pieces of helpful advice. 


\section{The Root Polytope}

\subsection{Kauffman States and the Alexander Polynomial}

To begin with, we describe a way to obtain a bipartite graph $G$ from a knot diagram $D$. First, construct the universe of $D$ (in the sense of \cite{Kau}) and color it in a checkerboard fashion.
Let us call the two colors black and white.
Next, put a black vertex in a black region and a white vertex in a white region.
Finally, connect the two vertices by an edge $e^*$ if the two regions that have these vertices share an edge $e$ of $D$.
In short, we obtain a bipartite graph from a knot diagram by considering the dual graph of the universe.

\begin{figure}[H] \label{bip}
 \centering
 \includegraphics[width=5cm]{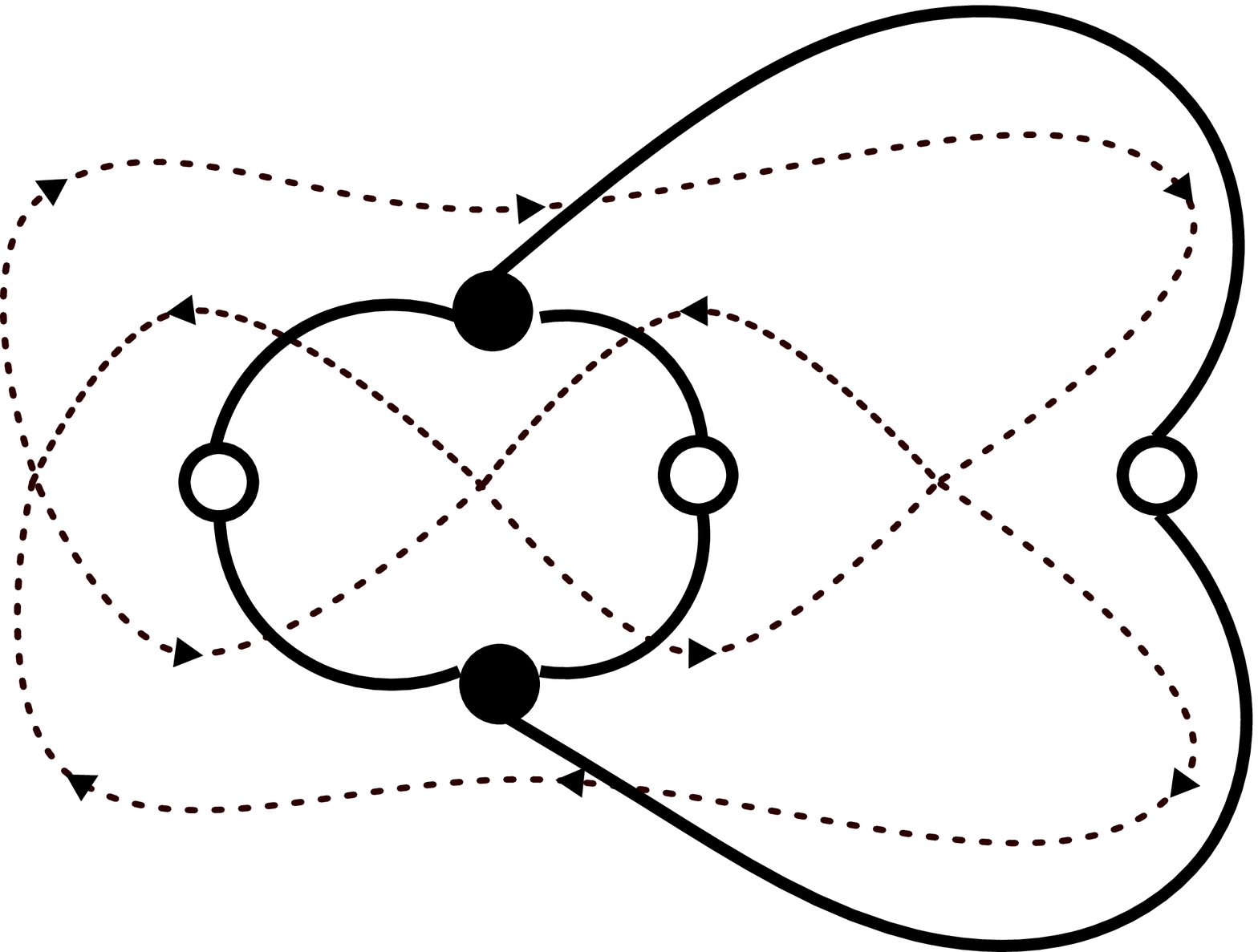}
 \hspace{5mm}
 \includegraphics[width=4cm]{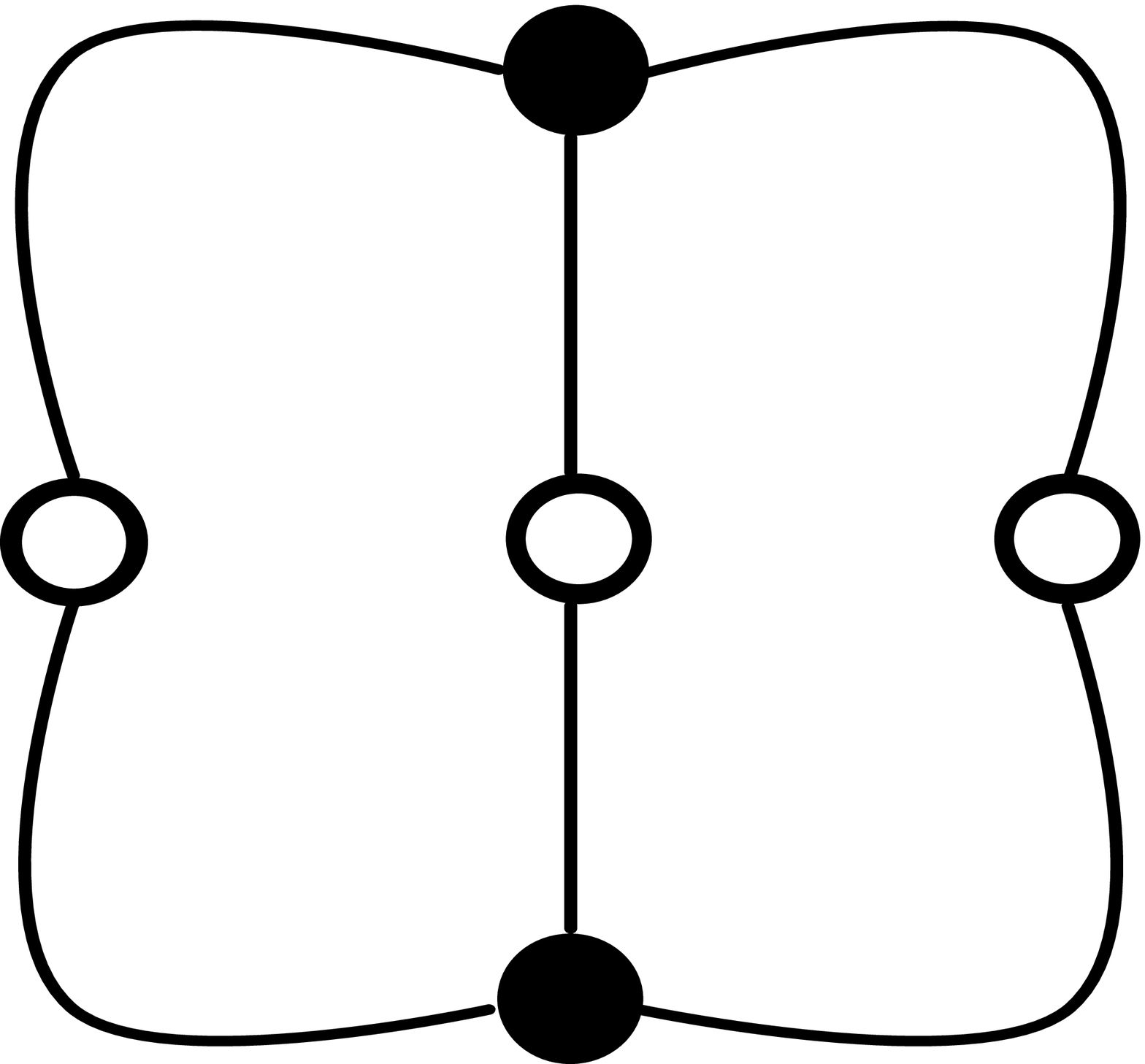}
 \caption{A universe (shown with dotted edges) and the corresponding bipartite graph.}
\end{figure}

We put an orientation on $D$ so that each edge of $D$ has a black vertex on its right side and a white one on the left side. Note that $D$ is balanced, i.e., the number of in-edges is equal to the number of out-edges at each vertex of $D$ (namely, both are $2$).

\begin{figure}[H] \label{Orientation}
 \centering
 \includegraphics[width=5cm]{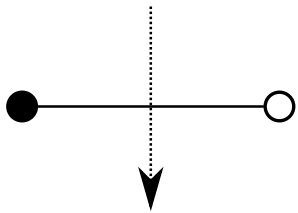}
 \caption{Orientation of the dual of a plane bipartite graph. The dotted line denotes an edge of $G^*$ and the solid line denotes an edge of $G$.}
\end{figure}

Provided a bipartite graph is given, following Postnikov \cite{P} we can construct a certain polytope in a Euclidean space. This polytope is called the root polytope.

\begin{defi}[Root polytope]
 Let $G$ be a bipartite graph with color classes $E$ and $V$. For $e \in E$ and $v \in V$, let $\mb{e}$ and $\mb{v}$ denote the corresponding standard generators of $\mathbb{R}^E \oplus \mathbb{R}^V$. 
 Define the root polytope of $G$ by
 \[ Q_G = \conv \{ \mb{e} + \mb{v} \mid ev \1 \mathrm{is} \1 \mathrm{an} \1 \mathrm{edge} \1 \mathrm{of} \1 G \}, \]
 where $\conv $ denotes the convex hull.
\end{defi}

If $G$ is connected, the dimension of $Q_G$ is $|E| + |V| -2$.

\begin{exam}
 Let $G$ be the bipartite graph in the right side of Figure~$1$. Since there are $6$ edges in $G$, we can plot $6$ vertices in five-dimensional Euclidean space. Taking the convex hull of these $6$ points, we obtain the root polytope
 shown in Figure~3. 
 
 \begin{figure}[h] \label{polytope}
  \centering
  \includegraphics[width=5cm]{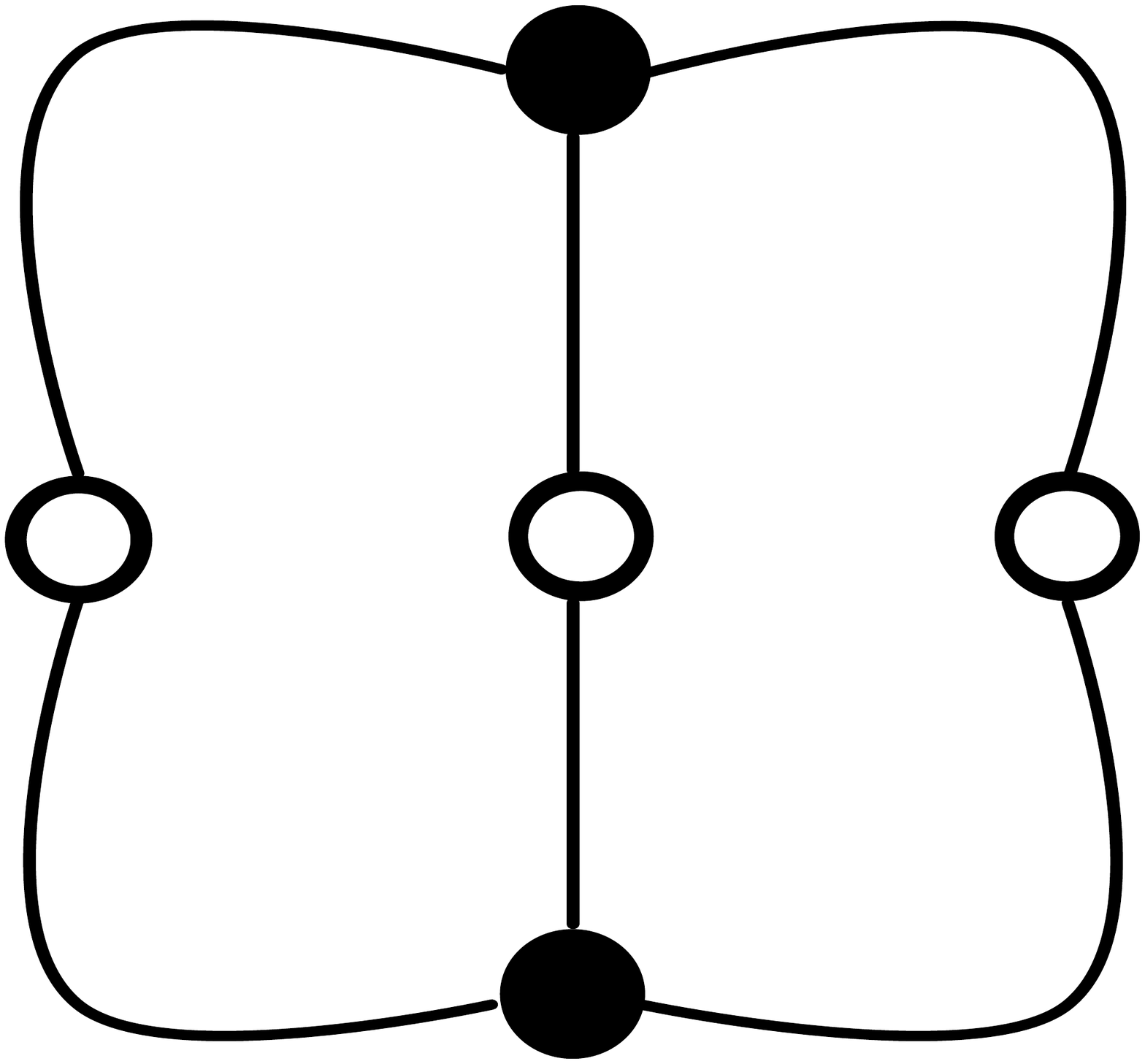}
  \includegraphics[width=5cm]{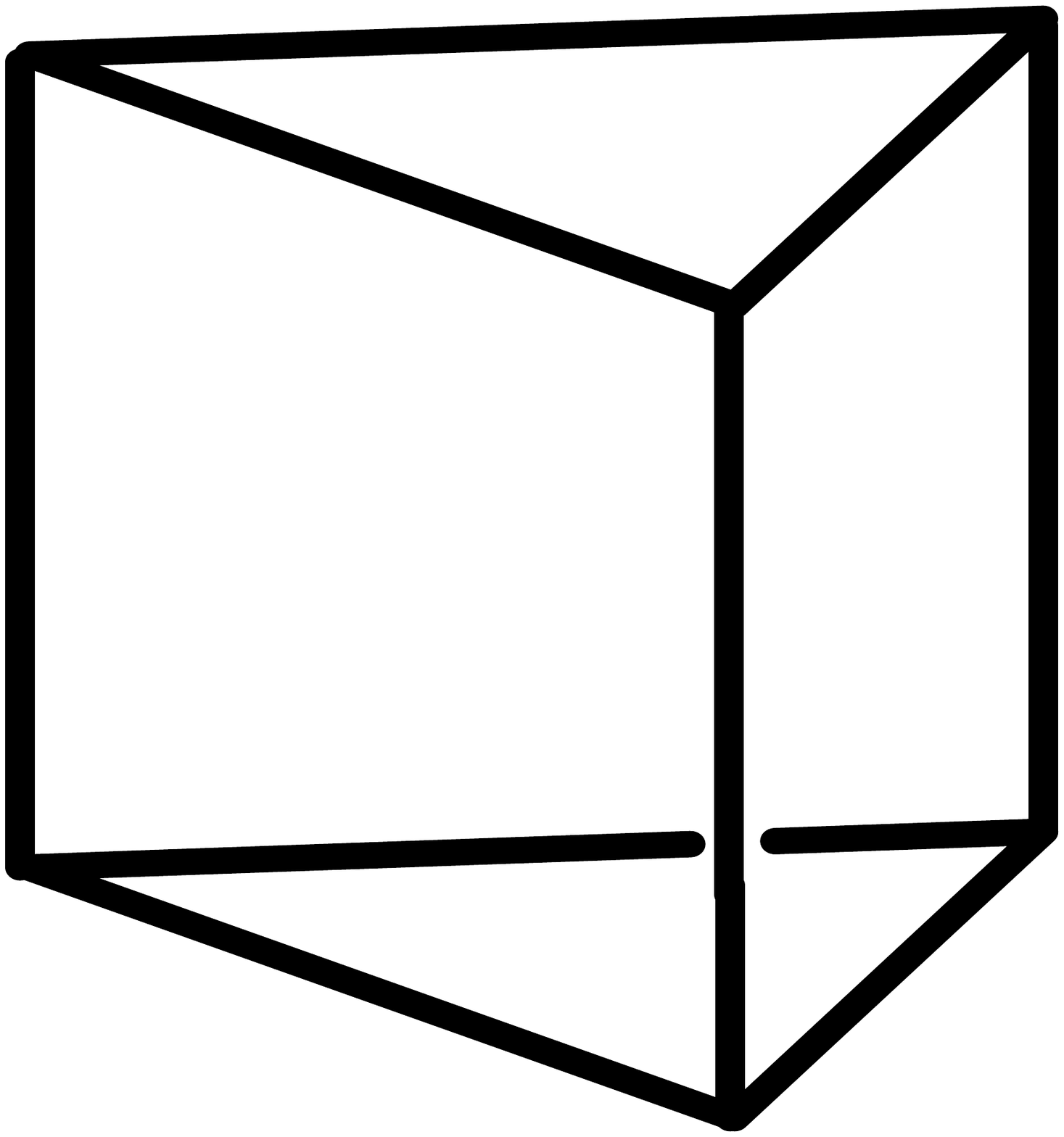}
  \caption{Bipartite graph (left) and its root polytope (right).}
 \end{figure}
 
\end{exam}

Let us introduce spanning arborescences to state our result. 

\begin{defi}[spanning arborescence]
Let $D$ be a directed graph.
Fix a vertex $r$ of $D$ and call it the root. A spanning tree in $D$ is called a spanning arborescence rooted at $r$ if the unique path in the tree from $r$ to any other vertex is oriented toward the root. The number of spanning arborescences is called the arborescence number of $D$ with respect to $r$.
\end{defi}

According to a combinatorial result \cite{S}, for any Eulerian directed graph $D$, the number of spanning arborescences in $D$ rooted at $r$, denoted by $\tau (D,r)$,  is given by the formula 
\[ \varepsilon (D,e) = \tau (D,r) \prod_{u\in R} \big( \mathrm{outdeg} (u) -1 \big)!, \]
where $\varepsilon (D, e)$ is the number of Eulerian tours in $D$ (starting from a fixed edge $e$) and $R$ is the set of vertices.

In our setting, the out-degrees on the right hand side are all equal to $2$.
Hence the formula says that the number of spanning arborescences is equal to the number of Eulerian tours.

This number equals the number of Kauffman states of $D$ with fixed stars. Here, a Kauffman state is an assignment of one marker per vertex (see Figure~4) so that each region in the universe receives no more than one marker.
As stated in Remark \ref{markers} below, we put stars in two adjacent regions and require the starred regions to be free of markers. These stars should be fixed throughout the argument.

\begin{figure}[H]
 \centering
 \includegraphics[width=3cm]{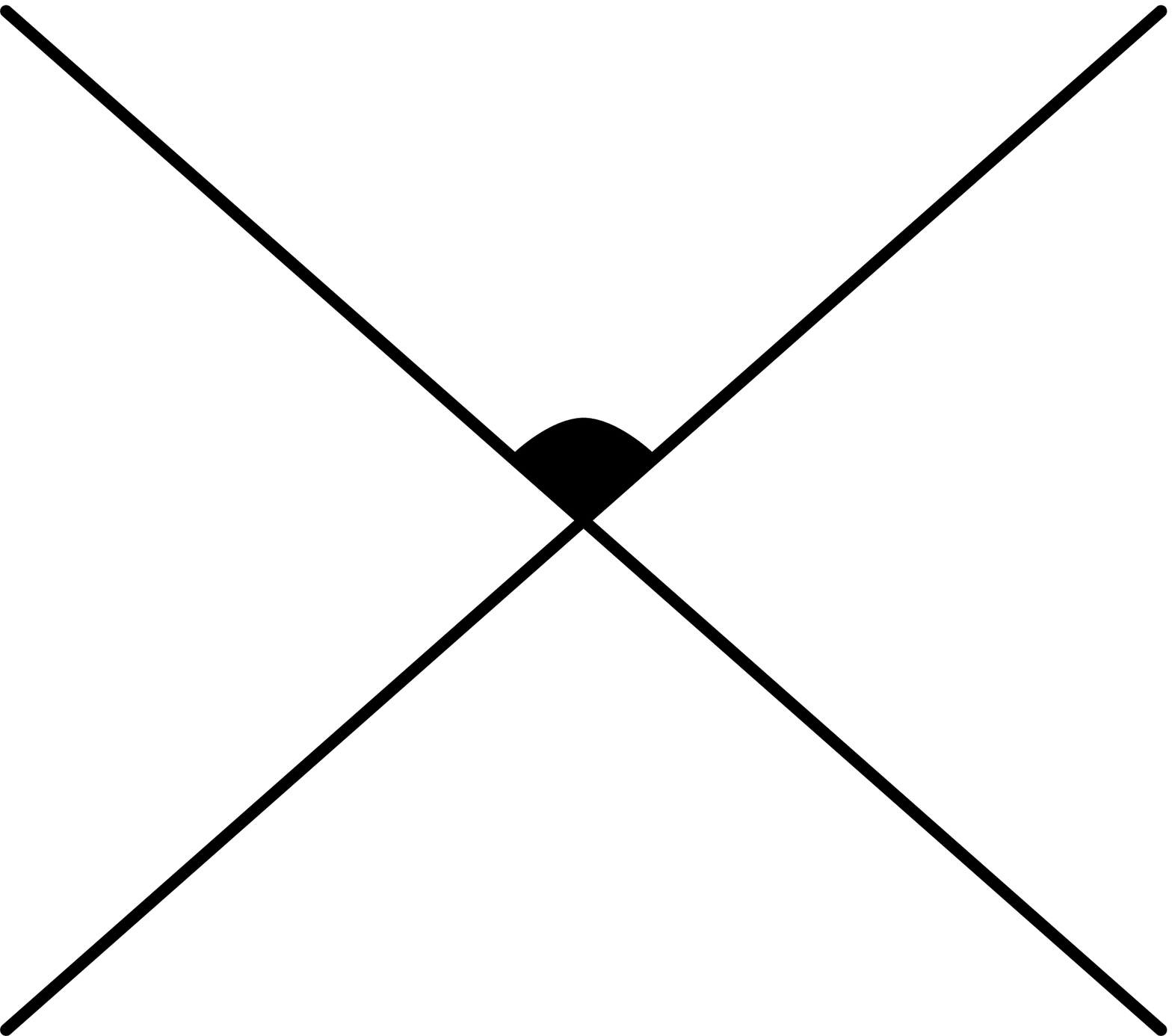}
 \caption{Marker.}
\end{figure}

\begin{rem} \label{markers}
 If the number of vertices is $n$, then the number of regions is $n+2$ by an application of Euler's formula. Since the number of the regions exceeds the number of vertices by $2$, 
 there are two regions that do not have a marker.
\end{rem}

Moreover, if $D$ is alternating, the number of Kauffman states coincides with the determinant $\det (L)$. We briefly recall the argument below.

\begin{figure}[H]
 \centering
 \labellist
  \pinlabel $U_k$ at 210 300
  \pinlabel $B_k$ at 100 200
  \pinlabel $W_k$ at 320 200
  \pinlabel $D_k$ at 210 100
 \endlabellist
 \includegraphics[width=4cm]{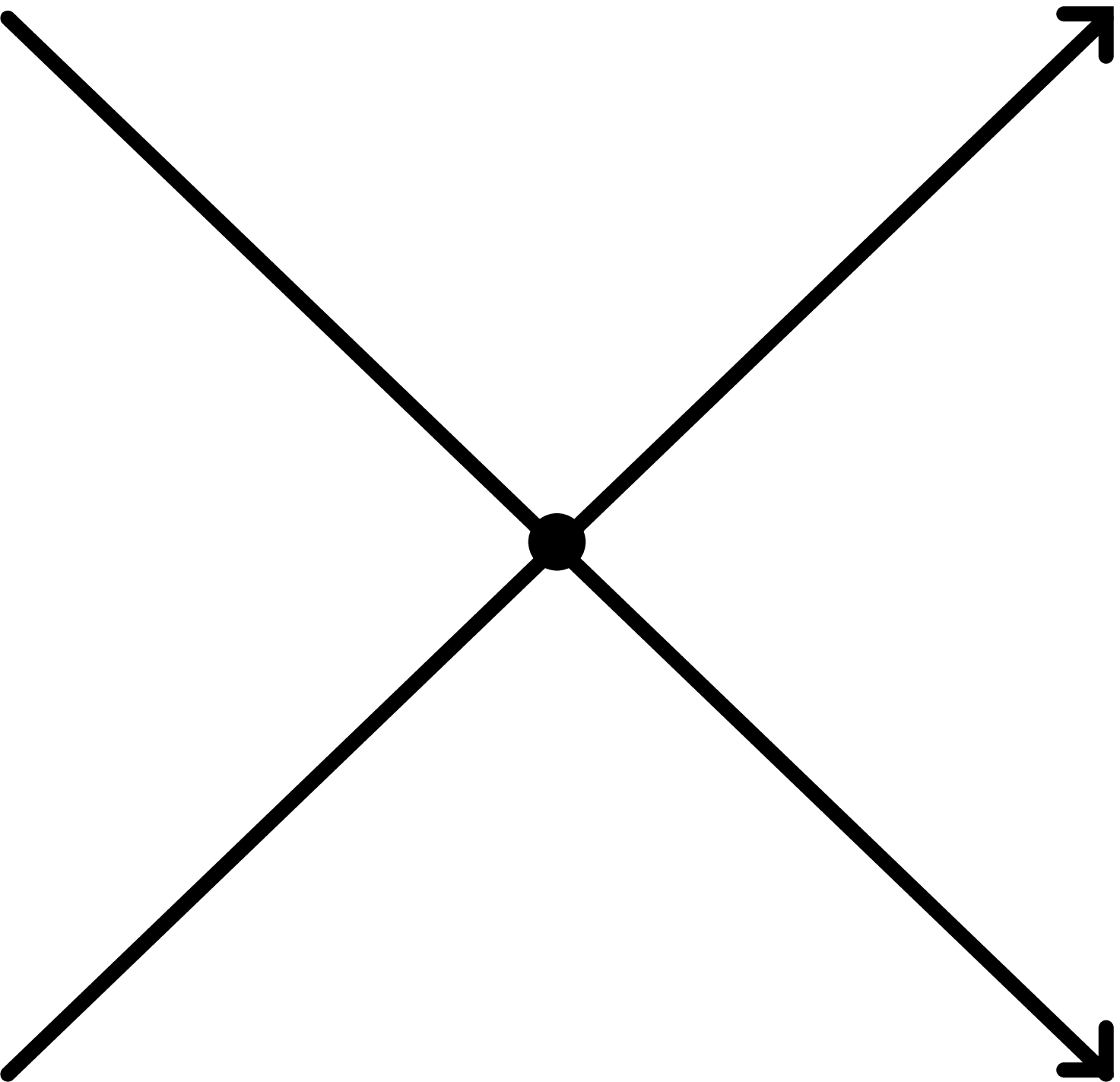}
 \caption{Labels at the $k^{\mathrm{th}}$ vertex.}
 \label{labels}
\end{figure}

Let us remember Kauffman's formula for the Alexander polynomial. To begin with, suppose that the universe has $n$ vertices. We give labels to the universe at the $k^{\mathrm{th}}$ vertex as shown in Figure \ref{labels}.

\begin{defi}
Let $K$ be a knot with labelled universe $D$ and $S$ be a Kauffman state of $D$. Define an inner product between $D$ and $S$ by
\[ \langle D \mid S \rangle = (-1)^{b(S)} V_1(S) V_2(S) \dots V_n(S), \]
where $b(S)$ denotes the number of black holes and $V_i(S)$ denotes the label touched by the marker at $i^{th}$ vertex.
Here, black hole is a marker that touches a region $B_k$ in Figure~5.
\end{defi}

We regard the inner product $\langle K\mid S \rangle$ as an element of the polynomial ring whose generators are the labels of $D$. Then we define a state polynomial for a labeled universe $D$.

\begin{defi} Let
 \[ \langle D \mid \mathfrak{S} \rangle = \sum_{S \in \mathfrak{S}} \langle K\mid S \rangle , \]
 where $\mathfrak{S}$ denotes the set of all Kauffman states of $D$.
\end{defi}

\begin{figure}[H]
 \centering
 \labellist
 \pinlabel $1$ at 400 360
 \pinlabel $-1$ at 340 310
 \pinlabel $t$ at 400 270
 \pinlabel $-t$ at 450 310
 \pinlabel $-1$ at 400 150
 \pinlabel $t$ at 350 90
 \pinlabel $1$ at 460 90
 \pinlabel $-t$ at 400 50
 \endlabellist
 \includegraphics[width=7cm]{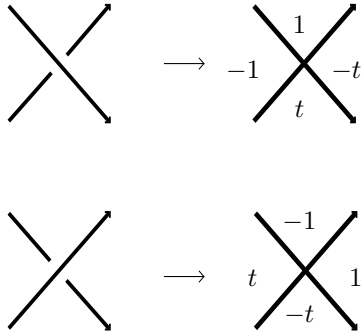}
 \caption{Labelings.}
 \label{labelings}
\end{figure}

\begin{theo}[\cite{Kau}]
 With an appropriate choice of labels, shown in Figure \ref{labelings},  the Kauffman state sum is the Alexander polynomial $\Delta_{K}(t)$.
\end{theo} 

To compute the determinant, we substitute $-1$ for $t$ in $\Delta_K$ and take the absolute value. 

\begin{figure}[H]
 \centering
 \labellist
 \pinlabel $\longrightarrow$ at 640 130
 \endlabellist
 \includegraphics[width=5cm]{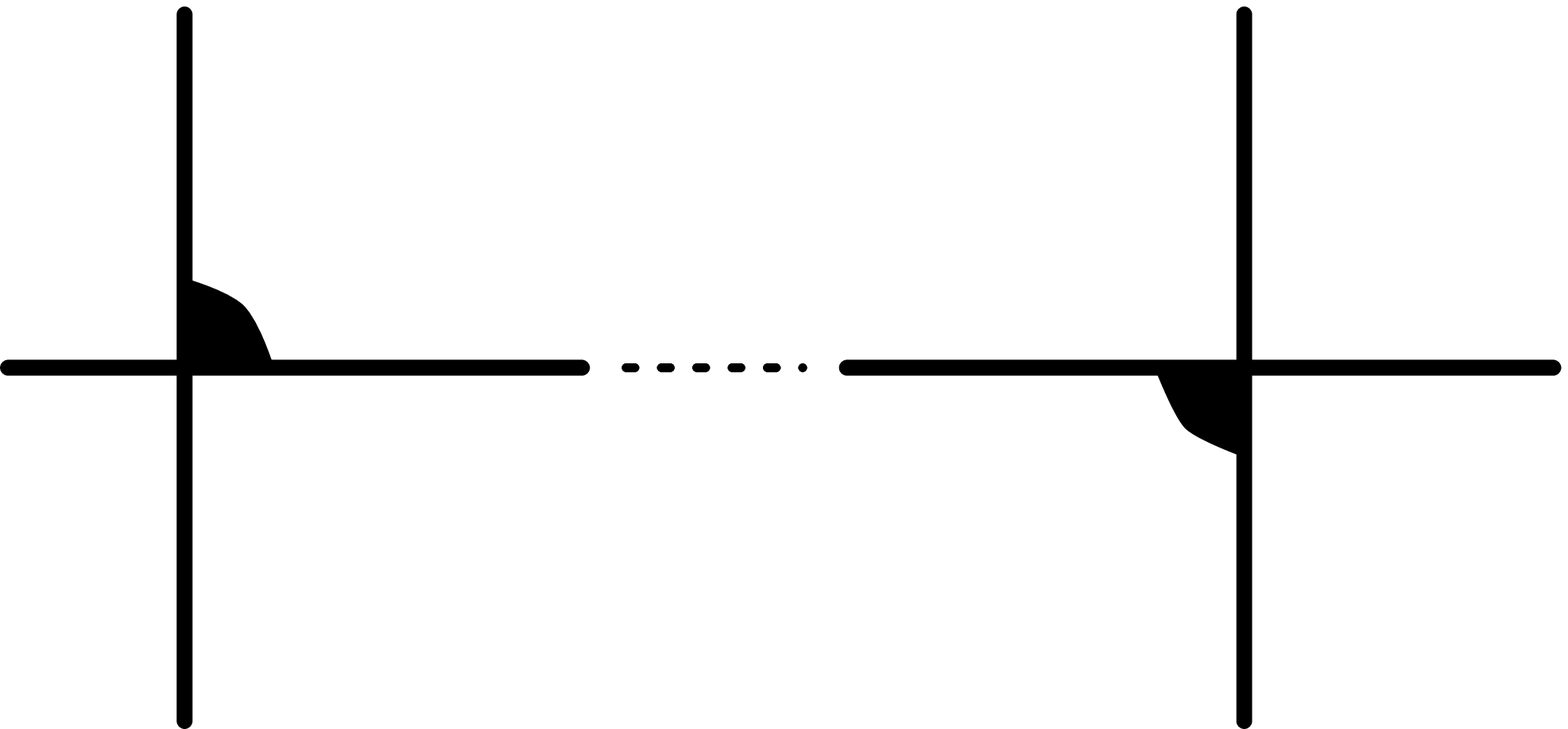}
 \hspace{1cm}
 \includegraphics[width=5cm]{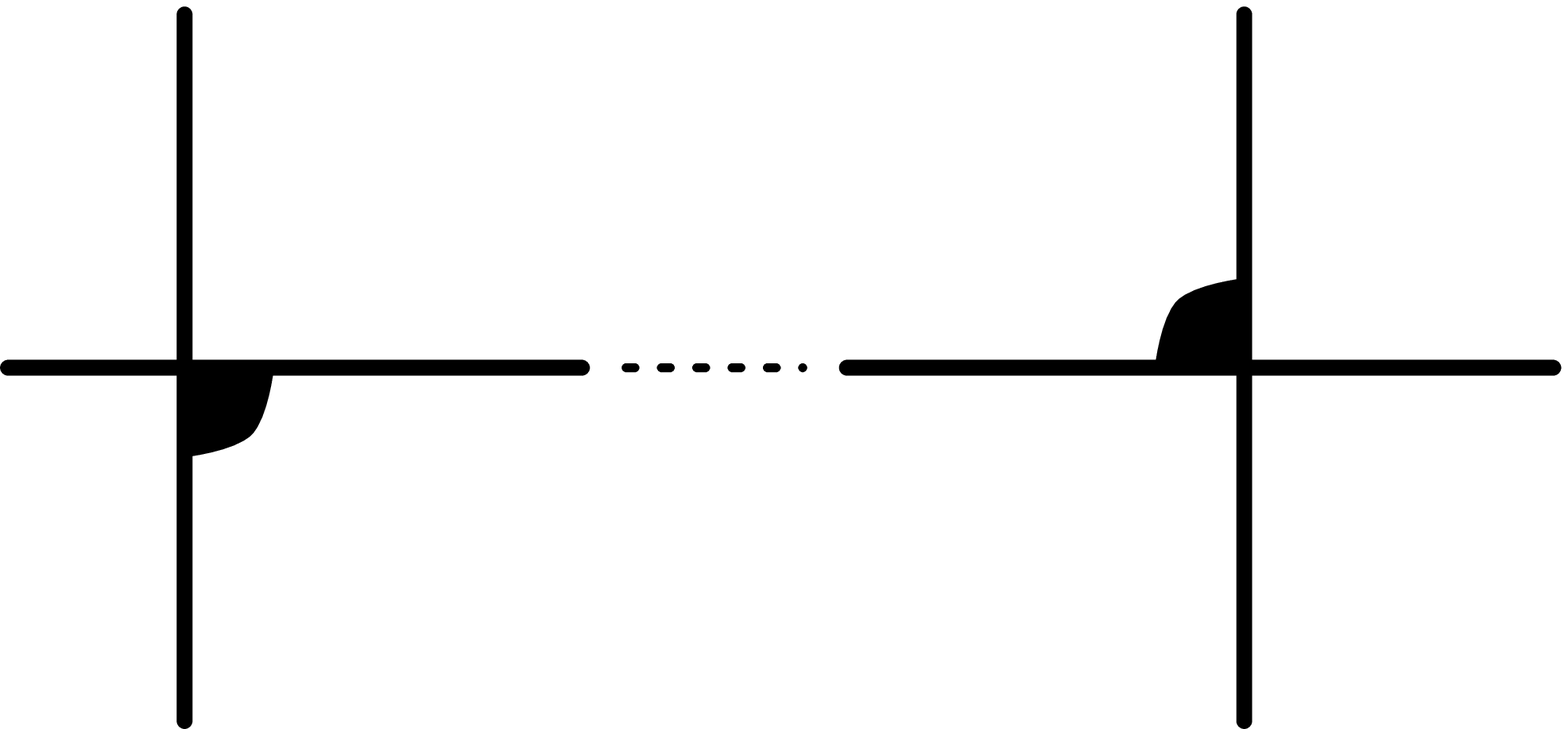}
 \caption{Clock move.}
\end{figure}

Kauffman's clock theorem \cite{Kau} says that any state of $D$ can be reached from the so-called clocked state by a sequence of clockwise moves.
By the definition of clockwise (counterclockwise) moves, we easily see that each clockwise (counterclockwise) move changes the sign $(-1)^{b(S)}$ of a state by $-1$ because the number of black holes changes by 1. Moreover, in case the given knot is alternating, the ratio of the contributions of the labels touched by the two markers to $\Delta _K(-1)$ is also $-1$. For details, see \cite{Kau}.
Putting them together, we can say that for an alternating link,  the inner product of a state with the Alexander labeling at $t=-1$ is not changed by clock moves.

\begin{figure}[H]
 \centering
 \labellist
  \pinlabel $\rightarrow$ at 730 155
  \pinlabel $-1$ at 878 230
  \pinlabel $1$ at 950 230
  \pinlabel $1$ at 1120 230 
  \pinlabel $1$ at 1190 230
  \pinlabel $-1$ at 878 170
  \pinlabel $1$ at 950 170
  \pinlabel $-1$ at 1120 170
  \pinlabel $-1$ at 1200 170
 \endlabellist
 \includegraphics[width=6cm]{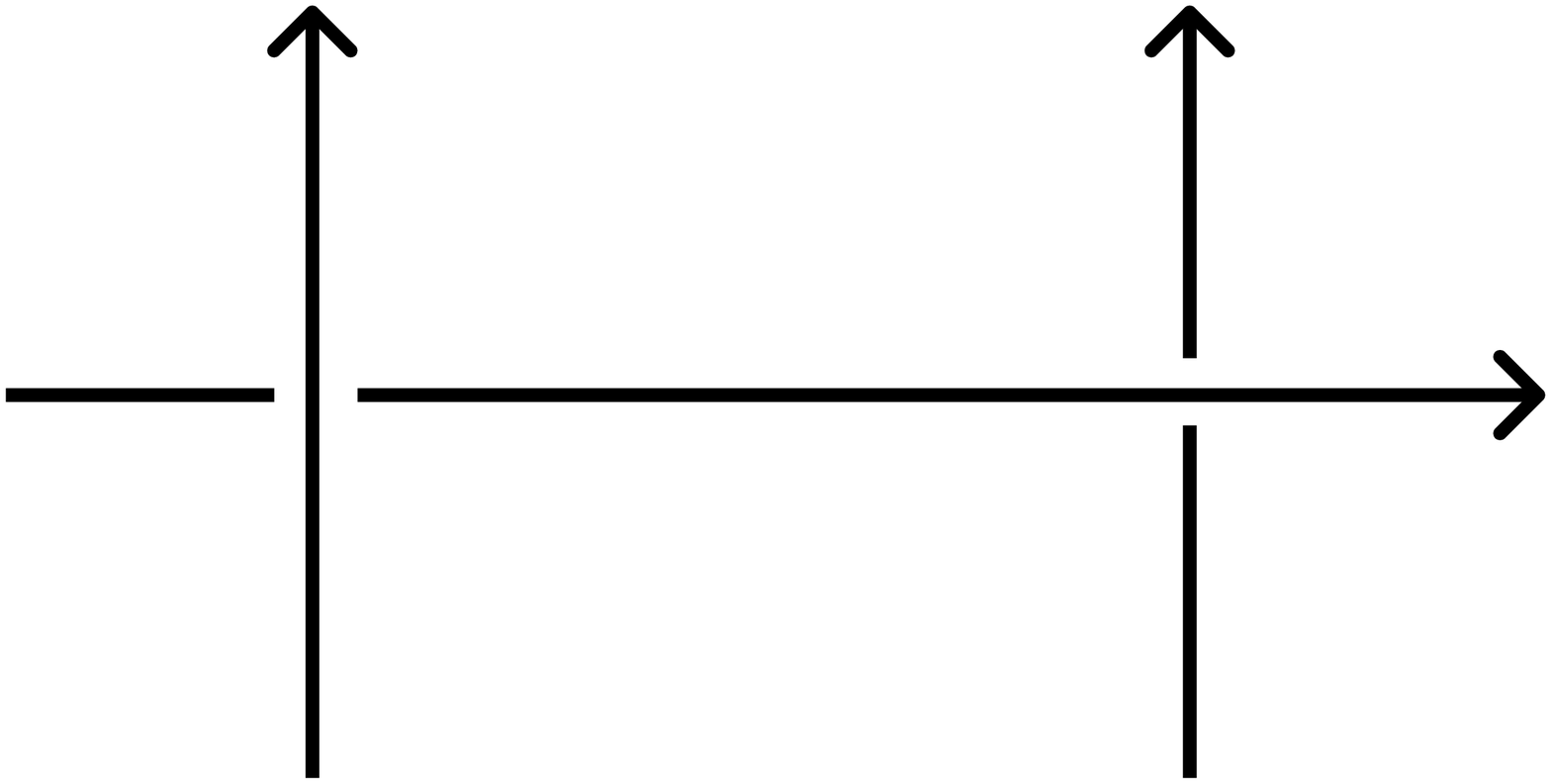}
 \hspace{1cm}
 \includegraphics[width=5cm]{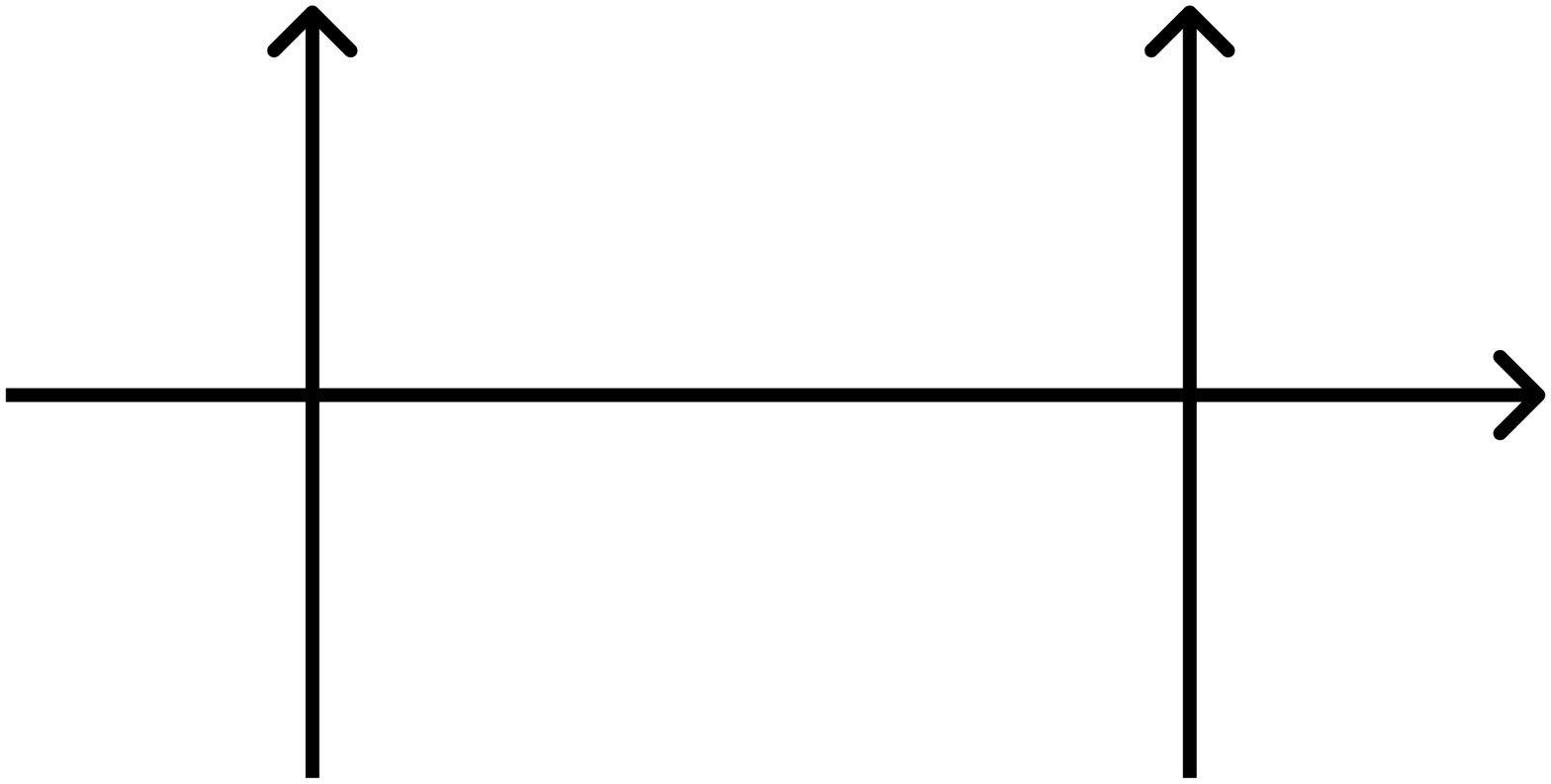}
 \caption{An example of labels in the alternating case.}
\end{figure}

\vspace{-5mm}

Since the clock theorem guarantees that each state of a given universe can be obtained from the clocked state, 
each state of the universe has the same inner product, namely $+1$ (or $-1$).

Now we can conclude that the number of Kauffman states is equal to plus/minus the value of the Alexander polynomial at $-1$, which is nothing but 
the determinant of the knot.

\subsection{Root polytope}

To connect the determinant and the bipartite graph $G=D^*$, we return to root polytopes.
Let us recall that a triangulation of the root polytope is a collection of maximal simplices in $Q_G$ so that their union is $Q_G$ and the intersection of any two simplices is their common face.

\begin{figure}[H]
 \centering
 \labellist
 \small
  \pinlabel $\mb{e_2+v_1}$ at 10 400
  \pinlabel $\mb{e_0+v_1}$ at 240 270
  \pinlabel $\mb{e_1+v_1}$ at 300 380
  \pinlabel $\mb{e_2+v_0}$ at 10 150
  \pinlabel $\mb{e_0+v_0}$ at 140 0
  \pinlabel $\mb{e_1+v_0}$ at 300 100
 \endlabellist
 \includegraphics[width=4cm]{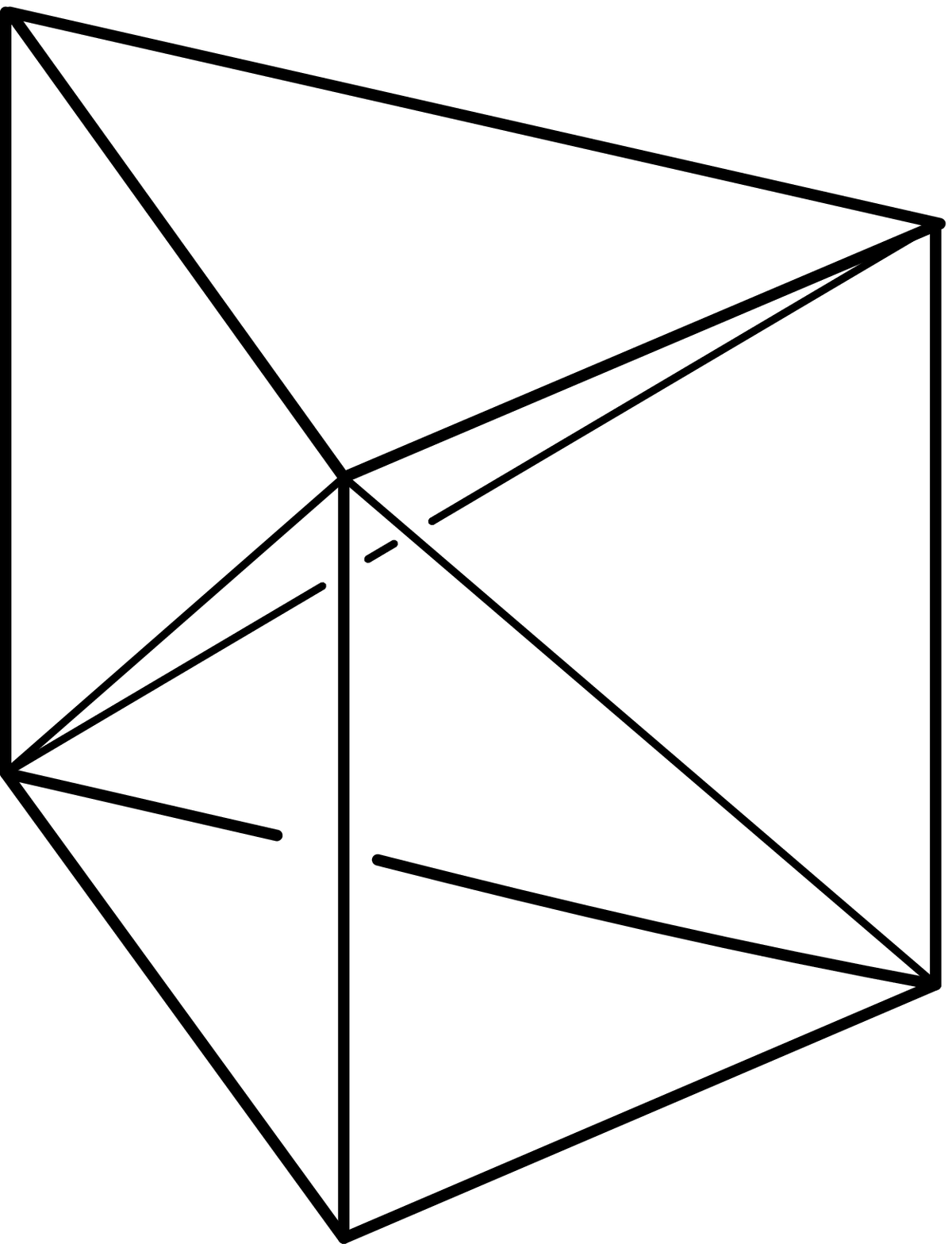}
 
 \vspace{4mm}
 
 \caption{Triangulation of the root polytope.}
\end{figure}

In \cite{KM2}, a way of triangulating the root polytope by means of spanning arborescences of $D$ is introduced.
For example, let us have a look at the spanning arborescence in Figure~$10$. 
This spanning arborescence has $2$ edges. Taking the dual of this arborescence, we obtain a spanning tree of $4$ edges. Then the corresponding polytope has $4$ vertices as shown in the figure.
For the reason that there are other $2$ spanning arborescences, this prism has $2$ more corresponding simplices to be triangulated.

Arborescences give us one triangulation. But the number of maximal simplicesin each triangulation of $Q_G$ is the same for the simple reason that all eligible simplices have the same volume \cite{P}.
Hence we obtain

\begin{theo}
Let $G$ be a plane bipartite graph with color classes $E$ and $V$ and $D=G^*$. Then the number of simplices in each triangulation of the root polytope $Q_G$ is the number of spanning arborescences 
in $G^*$.
\end{theo}

Since we have already seen that the number of spanning arborescences is equal to the determinant of the alternating link, 
we have the following theorem.

\begin{theo}
 The number of simplices to triangulate $Q_G$ equals the determinant of the alternating link $K$.
\end{theo}

\begin{figure}[H]
 \centering
 \labellist
 \small
  \pinlabel \textcolor{white}{$v_0$} at 310 455 
  \pinlabel \textcolor{white}{$v_1$} at 320 60
  \pinlabel $e_0$ at 130 250
  \pinlabel $e_1$ at 310 250
  \pinlabel $e_2$ at 480 250
  \pinlabel $\mathbf{e_2+v_1}$ at 870 500
  \pinlabel $\mathbf{e_1+v_1}$ at 1800 400
  \pinlabel $\mathbf{e_0+v_1}$ at 1100 250
  \pinlabel $\mathbf{e_0+v_0}$ at 870 40
 \endlabellist
 \includegraphics[width=4cm]{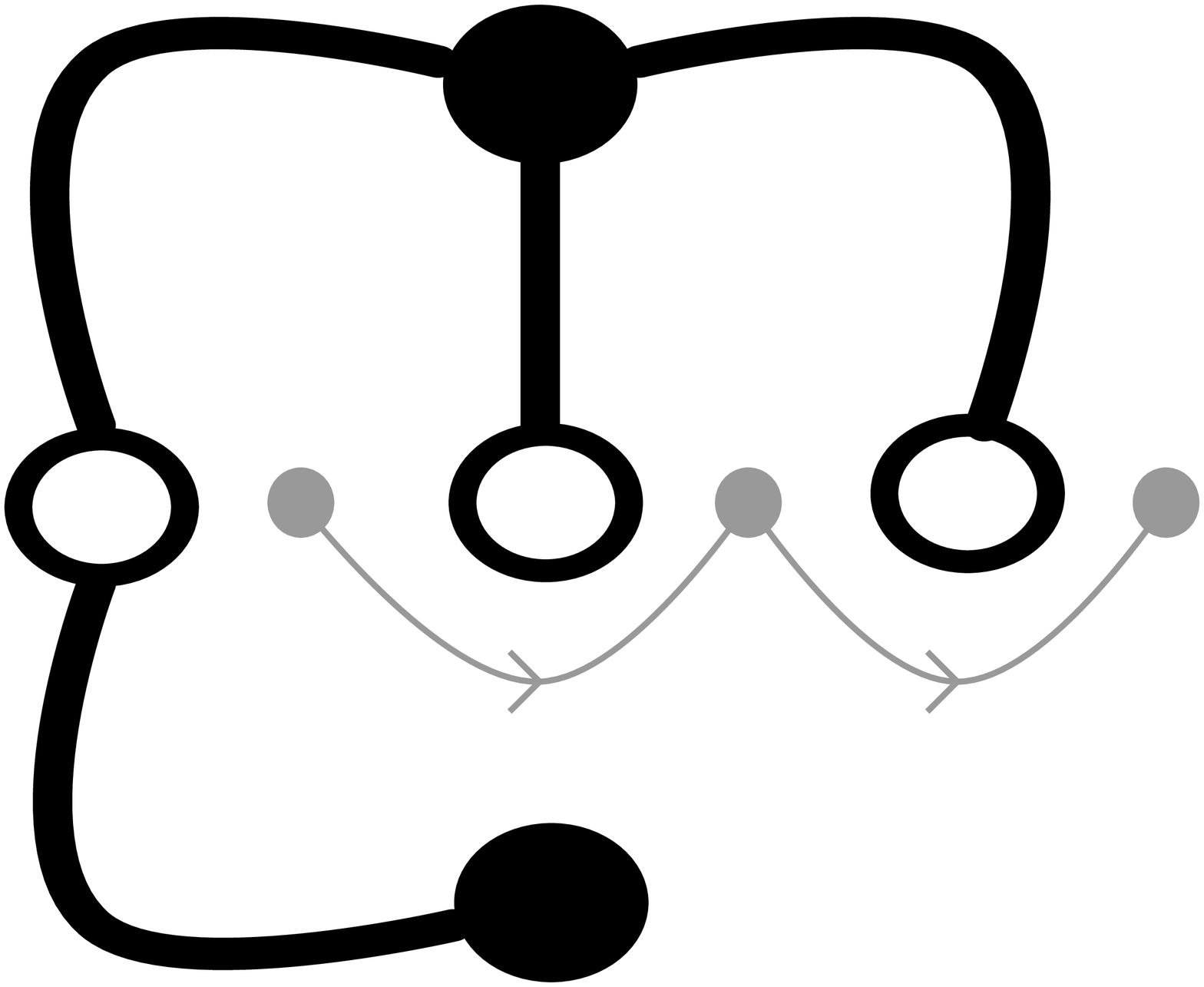}
 \hspace{15mm}
 \includegraphics[width=5cm]{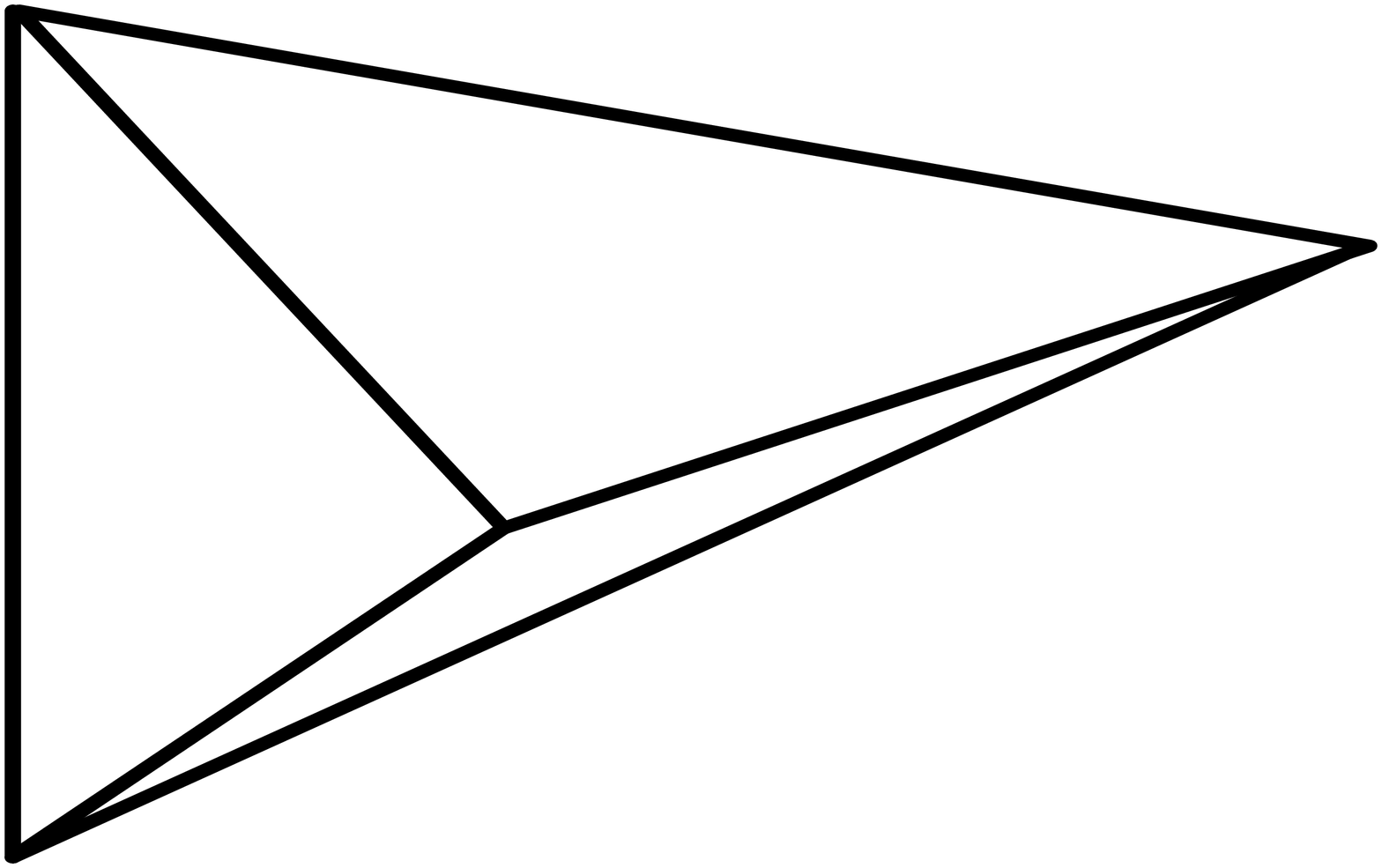}
 \caption{Spanning tree and corresponding simplex.}
\end{figure}

In \cite{K}, it is shown that the number of the spanning arborescences in $G^*$ is equal to the number of hypertrees in $G$. 
Putting all things together, we obtain that the number of hypertrees equals the determinant of the link $K$.
This completes the proof of Theorem 1.1.

Finally, we end with an open problem.

\begin{prob}
 Can the Alexander polynomial $\Delta_K$ or the Jones polynomial $J_K$ be expressed in terms of hypertrees?
\end{prob}




\begin{thebibliography}{99}
 \bibitem{B} A. Barvinok, Integer points in polyhedra, European Mathematical Society, 2008
 \bibitem{K} T. K\'alm\'an, A version of Tutte's polynomial for hypergraphs, Advances in Mathematics 244 (2013) 823-873.
 \bibitem{KM2} T. K\'alm\'an and Hitoshi Murakami, Root polytopes, parking functions, and the HOMFLY polynomial, to appear in Quantum Topology.
 \bibitem{Kau} L. Kauffman, Formal knot theory, Mathematical Notes 30 Princeton University Press 1983.
 \bibitem{L} R. Lickorish, An introduction to knot theory, Graduate Texts in Mathematics, vol. 175, Springer-Verlag, New York, 1997.
 \bibitem{P} A. Postnikov, Permutohedra, associahedra, and beyond,  International Mathematics Research Notices.
 \bibitem{S} R. Stanley, Enumerative Combinatorics 2, Cambridge Studies in Advanced Mathematics 62.
\end{thebibliography}
\end{document}